\documentclass{amsart}

\usepackage{amsfonts,amssymb,verbatim,amsmath,amsthm,latexsym,textcomp,amscd}
\usepackage{latexsym,amsfonts,amssymb,epsfig,verbatim}
\usepackage{amsmath,amsthm,amssymb,latexsym,graphics,textcomp}
\usepackage{graphicx}
\usepackage{color}
\usepackage{url}

\begin{document}

\newtheorem{theorem}{Theorem}[section]
\newtheorem{prop}[theorem]{Proposition}
\newtheorem{lemma}[theorem]{Lemma}
\newtheorem{cor}[theorem]{Corollary}
\newtheorem{defn}[theorem]{Definition}
\newtheorem{conj}[theorem]{Conjecture}
\newtheorem{claim}[theorem]{Claim}
\newtheorem{example}[theorem]{Example}
\newtheorem{rem}[theorem]{Remark}
\newtheorem{rmk}[theorem]{Remark}
\newtheorem{obs}[theorem]{Observation}
\newtheorem{qn}[theorem]{Question}

\newcommand{\map}{\rightarrow}
\newcommand{\C}{\mathcal C}
\newcommand\AAA{{\mathcal A}}
\newcommand\BB{{\mathcal B}}
\newcommand\DD{{\mathcal D}}
\newcommand\EE{{\mathcal E}}
\newcommand\FF{{\mathcal F}}
\newcommand\GG{{\mathcal G}}
\newcommand\HH{{\mathcal H}}
\newcommand\I{{\stackrel{\rightarrow}{i}}}
\newcommand\J{{\stackrel{\rightarrow}{j}}}
\newcommand\K{{\stackrel{\rightarrow}{k}}}
\newcommand\LL{{\mathcal L}}
\newcommand\MM{{\mathcal M}}
\newcommand\NN{{\mathbb N}}
\newcommand\OO{{\mathcal O}}
\newcommand\PP{{\mathcal P}}
\newcommand\QQ{{\mathcal Q}}
\newcommand\RR{{\mathcal R}}
\newcommand\SSS{{\mathcal S}}
\newcommand\TT{{\mathcal T}}
\newcommand\UU{{\mathcal U}}
\newcommand\VV{{\mathcal V}}
\newcommand\WW{{\mathcal W}}
\newcommand\XX{{\mathcal X}}
\newcommand\YY{{\mathcal Y}}
\newcommand\ZZ{{\mathbb Z}}
\newcommand\hhat{\widehat}
\newcommand\vfn{\stackrel{\A}{r}(t)}
\newcommand\dervf{\frac{d\stackrel{\A}{r}}{dt}}
\newcommand\der{\frac{d}{dt}}
\newcommand\vfncomp{f(t)\I+g(t)\J+h(t)\K}
\newcommand\ds{\sqrt{f^{'}(t)^2+g^{'}(t)^2+h^{'}(t)^2}dt}
\newcommand\rvec{\stackrel{\A}{r}}
\newcommand\velo{\frac{d\stackrel{\A}{r}}{dt}}
\newcommand\speed{|\velo|}
\newcommand\velpri{\rvec \,^{'}}
\newcommand{\RED}{\textcolor{red}}

\title[Corrigendum]{Corrigendum to ``Graphs of hyperbolic groups and a limit set intersection theorem"}

\author{Pranab Sardar}
\address{Indian Institute of Science Education and Research Mohali}

\date{\today}

\maketitle

\section{Introduction}
The purpose of this note is to point out a mistake in the proof of the Proposition $4.9$ of \cite{ps-limitset} which subsequently weakens
the main results of \cite{ps-limitset}- Theorem $4.1$ and Corollary $4.10$. We state and prove the weaker statements here. 
See Theorem \ref{main-thm-v2} and Corollary \ref{final-cor-v2}. Finally we provide an example where Corollary $4.10$ of \cite{ps-limitset}
fails. 
\begin{defn}
$(1)$ Suppose $X$ is a (proper) hyperbolic metric space and $Y\subset X$. A point $\xi\in \Lambda(Y)\subset \partial X$ is a called
a {\em conical limit point} of $Y$ if for any geodesic ray $\alpha$ in $X$ such that $\alpha(\infty)=\xi$ there is a constant $D>0$
such that there is a sequence $\{y_n\}$ contained in $N_D(\alpha)\cap Y$ converging to $\xi$.

$(2)$ Suppose $H$ is a group acting on a hyperbolic metric space $X$ by isometries.
  A point $\xi\in \Lambda(H)$ is called a {\em conical limit point} of $H$ if for any (equivalently {\em some}) $x_0\in X$, 
$\xi$ is a conical limit point of the orbit $H.x_0$.

$(3)$ The set of all conical limit points of $H$, denoted by $\Lambda_c(H)$, will be called the {\em conical limit set} of $H$.
\end{defn}

{\bf Remark:} The second part of the above definition makes sense for any infinite subset of the group $H$ too. We are interested in the
case where $X$ is a Cayley graph of a Gromov hyperbolic group $G$ on which $G$ has the natural action and $H$ is a subgroup of $G$
or more generally a subset of $G$. The following lemma then follows as that of Lemma $1.2$ of \cite{ps-limitset}.

\begin{lemma}\label{int-lemma-v2}
Suppose $G$ is a hyperbolic group and $H$ is any {\em subset} of $G$. Then for all $x\in G$ we have

$(1)$ $\Lambda_c(xH)=\Lambda_c(xHx^{-1})$.

$(2)$ $\Lambda_c(xH)= x\Lambda_c(H)$.
\end{lemma}

\begin{defn}({\bf Limit set intersection property}, \cite{ps-limitset})
Suppose $G$ is a Gromov hyperbolic group.
Let $\mathcal S$ be any collection of subgroups of  $G$. We say that $\mathcal S$
has the {\em limit set intersection property} if for all $H, K\in \mathcal S$ we have
$\Lambda(H)\cap \Lambda(K)=\Lambda(H\cap K)$.
\end{defn}

\begin{defn}({\bf Conical limit set intersection property})
Suppose $G$ is a Gromov hyperbolic group.
Let $\mathcal S$ be any collection of subgroups of  $G$. We say that $\mathcal S$
has the {\em conical limit set intersection property} if for all $H, K\in \mathcal S$ we have
$\Lambda_c(H)\cap \Lambda_c(K)=\Lambda_c(H\cap K)$.
\end{defn}


\section{The main result}

For the rest of this note we shall assume all the notation and hypotheses of the section $4$ of \cite{ps-limitset}.
In particular we assume that $G$ is a hyperbolic group which admits a graph of groups decomposition
$(\mathcal G, Y)$ with the QI embedded condition where all the vertex and edge groups are hyperbolic and
$\mathcal T$ is the Bass-Serre tree for this graph of groups. We denote by $X$ the space quasi-isometric to $G$ as
constructed in section $3$ of \cite{ps-limitset} from the graph of groups decomposition of $G$ and 
$\Theta$ will denote an orbit map $G\map X$. As noted after Corollary $3.6$ of \cite{ps-limitset} 
$\Theta$ naturally induces to a (uniform) quasi-isometry from $gG_v$ (or $gG_vg^{-1}$)
to the corresponding vertex space $X_{\tilde{v}}$ where $\tilde{v}=gG_v$
for all $v\in V(Y)$ and $g\in G$. These maps will be loosely referred to as the {\em restrictions of $\Theta$}. We note that
under this map a conical limit point of $gG_v$ in $\partial G$ will be mapped to a conical limit point of $X_{\tilde{v}}$
in $\partial X$ since $\Theta$ is a quasi-isometry.

We prove the following weaker alternative for the
Theorem $4.1$ of \cite{ps-limitset}; namely that the vertex groups $\{G_v:v\in V(\TT)\}$ satisfy the conical limit set
intersection property:

\begin{theorem}\label{main-thm-v2}
Suppose a hyperbolic group $G$ admits a decomposition into a graph of hyperbolic groups $(\mathcal G,Y)$ with
quasi-isometrically embedded condition and suppose $\mathcal T$ is the corresponding Bass-Serre tree.
Then for all $w_1, w_2\in V(\mathcal T)$ we have 
$\Lambda_c(G_{w_1})\cap \Lambda_c(G_{w_2})= \Lambda_c(G_{w_1}\cap G_{w_2})$.
\end{theorem}

This result can be generalized to prove the corresponding analog of Corollary $4.10$ of \cite{ps-limitset} in the
same of way as Corollary $4.10$ was derived from Theorem $4.1$ in \cite{ps-limitset}.

\begin{cor}\label{final-cor-v2}
If $H_i\subset G_{w_i}$, $i=1,2$ are two quasiconvex subgroups then $\Lambda_c(H_1)\cap \Lambda_c(H_2)=\Lambda_c(H_1\cap H_2)$.
\end{cor} 


{\bf Comments on the proofs of Proposition $4.9$ and Theorem $4.1$ of \cite{ps-limitset}:}
\smallskip
Let $[v,w]$ denote the geodesic in $\TT$ joining $v,w$ for all $v,w\in V(\TT)$.
Suppose $\xi_v\in \partial X_v$ and $\xi_w \in \partial X_w$ map to the same point $\xi\in \partial X$ under the CT maps
$\partial {X}_v\rightarrow \partial {X}$ and $\partial {X}_w\rightarrow \partial {X}$ respectively. 

$(1)$ In the first sentence of the proof of Proposition $4.9$ we made the following tacit assumption which is wrong: 

{\em Suppose $\xi_v$ cannot be flowed to $X_w$ and $\xi_w$ cannot be flowed to $X_v$.
Suppose $v_1$ is the farthest vertex from $v$ on $[v,w]$ such that
$\xi_v$ can be flowed to $X_{v_1}$ and $w_1$ is the nearest vertex from $v$ on $[v,w]$ such that $\xi_w$ can be flowed to $X_{w_1}$.
Then we have $d_{\TT}(v,v_1)<d_{\TT}(v,w_1)$.}

It could as well happen that $d_{\TT}(v,v_1)\geq d_{\TT}(v,w_1)$ although $v_1\neq w$ and $w_1\neq v$. However, rest of the proof of
that proposition is correct with this tacit assumption. The current proof shows the following weaker statement:

\begin{prop}({\bf A correct alternative to Proposition $4.9$ of \cite{ps-limitset}})\label{old-prop-v2}
Suppose $\xi_v\in \partial X_v$ and $\xi_w \in \partial X_w$ map to the same point $\xi\in \partial X$ under the CT maps
$\partial {X}_v\rightarrow \partial {X}$ and $\partial {X}_w\rightarrow \partial {X}$ respectively. Then there is a vertex
$z\in [v,w]$ such that $\xi_v$ and $\xi_w$ both can be flowed to $X_z$.
\end{prop}

$(2)$ In \cite{ps-limitset} the proof of Theorem $4.1$ is dependent on Proposition $4.9$ which is mentioned 
in the last fourth sentence of the second paragraph of the proof. All that stated before in that proof are independent of this erroneous
proposition and the rest of the argument is also independently correct modulo this assumption. The current proof shows
the following weaker result:

\begin{theorem}({\bf A correct alternative to Theorem $4.1$ of \cite{ps-limitset}})\label{old-thm-v2}
Suppose $v,w\in V(\mathcal T)$ and $\xi\in \Lambda(G_v)\cap \Lambda(G_w)$.
Suppose $\xi_v\in \partial X_v$ and $\xi_w \in \partial X_w$ both map to $\partial \Theta(\xi)$
under the CT maps $\partial {X}_v\rightarrow \partial {X}$ and $\partial {X}_w\rightarrow \partial {X}$ respectively.
If either $\xi_v$ can be flowed to $X_w$ or $\xi_w$ can be flowed to $X_v$ then $\xi\in \Lambda(G_{v}\cap G_{w})$.
\end{theorem}


{\bf Proof of Theorem \ref{main-thm-v2}:} 

\smallskip

We first prove the following two propositions.

\begin{prop}\label{prop1}
Suppose $v,w\in V(\TT)$ and $\xi_v\in \partial X_v$ can be flowed to $\xi_w\in \partial X_w$.
Suppose $\xi\in \partial X$ is the image of both $\xi_v$ and $\xi_w$ under the CT maps $\partial X_v\map \partial X$
and $\partial X_w\map \partial X$ respectively.

If $\xi$ is a conical limit point of $X_v$ then $\xi$ is a conical limit point of $X_w$.
\end{prop}

Note that if $\xi_v\in \partial X_v$ can be flowed to $\xi_w\in \partial X_w$ then by Corollary $4.5$ of \cite{ps-limitset}
they are mapped to the same point of $\partial X$ under the CT maps $\partial X_v\map \partial X$
and $\partial X_w\map \partial X$.
\smallskip

{\em Proof of the proposition:} Suppose $\alpha_v\subset X_v$ and $\alpha_w\subset X_w$ are geodesic rays such that $\alpha_v(\infty)=\xi_v$
and $\alpha_w(\infty)=\xi_w$. By Lemma $4.4$ of \cite{ps-limitset} $Hd(\alpha_v,\alpha_w)<\infty$. Hence, it is enough to show
that if $\alpha$ is a geodesic ray in $X$ joining $\alpha(0)$ to $\xi$ then there is $D>0$ and an unbounded sequence of points
$\{x_n\}$ on $\alpha_v$ such that $d(x_n, \alpha)\leq D$. This is what we prove next.

Let $\lambda_n$ be the portion of $\alpha_v$ from $\alpha_v(0)$ to $\alpha_v(n)$. We know that
the ladder $B(\lambda_n)$ is uniformly quasiconvex in $X$ (see Theorem $4.7$, \cite{ps-limitset}). 
Hence there is a uniform quasigeodesic in $X$ contained in $B(\lambda_n)$
joining $\alpha_v(0)$ and $\alpha_v(n)$. Let $\gamma_n\subset B(\lambda_n)$ be such a quasigeodesic. 
Now since $x_n\map \xi$, $\gamma_n$ fellow travels with $\alpha$ for a long time for large $n$. 
Hence there is a uniform constant $D_1>0$ and a sequence of points $\{y_n\}$, $y_n\in \gamma_n$ for all $n$ 
such that $d(\alpha_v(0), y_n)\map \infty$, and the portion of $\gamma_n$ between $\alpha_v(0)$ and $y_n$ is contained
in $N_{D_1}(\alpha)$. Now since $\xi$ is a conical limit point there is a constant $D_2>0$ such that 
$N_{D_2}(\alpha)\cap X_v$ is an infinite set. Hence, there is an unbounded sequence of points $\{z_n\}$ on $\alpha$
such that $d(z_n, X_v)\leq D_2$ for all $n$. However, for all $n$ there is $n^{'}$ such that $d(z_n, \gamma_{n^{'}})\leq D_1$.
Let $z^{'}_n\in \gamma_{n^{'}}$ be such that $d(z_n, z^{'}_n)\leq D_1$. Then $d(z^{'}_n, X_v)\leq D_1+D_2$.
Let $z^{''}_n\in X_v$ be such that $d(z^{'}_n, z^{''}_n)\leq D_1+D_2$.
Now we use Mitra's projection $P:X\map B(\lambda_{n^{'}})$. Since $P(z^{'}_n)=z^{'}_n$, if we set $x_n=P(z^{''}_n)$ we
see that $d(x_n, z^{'}_n)$ is uniformly bounded. It follows that $d(x_n, z_n)$ uniformly bounded too. $\Box$

\begin{prop}\label{prop2}
Suppose $\tilde{v}=gG_v\in V(\TT)$ where $v\in V(Y)$ and $g\in G$. Suppose $\eta_1, \eta_2\in \partial X_{\tilde{v}}$
are mapped to the same point $\eta\in \partial X$ under the CT map
$\partial X_{\tilde{v}}\map \partial X$. If $\eta$ is a conical limit point of $X_{\tilde{v}}$ then $\eta_1=\eta_2$. 
\end{prop}

$Proof:$ Let $\xi\in \partial G$ and $\xi_i\in \partial gG_vg^{-1}$ be such that their images under the boundary maps induced
by $\Theta$ and its restriction to $gG_vg^{-1}$ are respectively $\eta$ and $\eta_i$, $i=1,2$. Then the proposition follows from
the second part of the Theorem A of \cite{jeon-etal}. $\Box$

{\bf Remark:} One can give an independent proof of Proposition \ref{prop2} in the line of the proof of Theorem $4.11$
of \cite{mitra-endlam} using Mitra's ladders which would work for any tree of hyperbolic metric spaces with QI embedding condition.
However, the proof would then be much longer.

\smallskip

{\em Proof of Theorem \ref{main-thm-v2}:} Clearly $\Lambda_c(G_{w_1}\cap G_{w_2})\subset \Lambda_c(G_{w_1})\cap \Lambda_c(G_{w_2})$.
Suppose $\xi\in \Lambda_c(G_{w_1})\cap \Lambda_c(G_{w_2})$. Let $\xi_i\in \partial G_{w_i}$, $i=1,2$ be such that they both map
to $\xi$ under the CT maps $\partial G_{w_i}\map \partial G$. Let $\eta\in \partial X$ and $\eta_i\in \partial X_{w_i}$ be the
images of $\xi,\xi_1, \xi_2$ respectively under the image of $\Theta$ and its restrictions to $G_{w_1}$ and $G_{w_2}$ respectively.
Clearly $\eta_i$ maps to $\eta$ under the CT maps $\partial X_{w_i}\map \partial X$, $i=1,2$.
Hence, by Proposition \ref{old-prop-v2} there is a vertex $w\in [w_1,w_2]$ such that both $\eta_i$ can be flowed to $X_w$.
Let the flowed images be $\eta^{'}_i\in \partial X_w$ respectively. Then $\eta^{'}_i$'s also map to $\eta$ under the CT map 
$\partial X_w\map \partial X$ by Corollary $4.5$ of \cite{ps-limitset}.
By Lemma \ref{int-lemma-v2} $\eta$ is a conical limit point of $X_{w_i}$, $i=1,2$. Then
by Proposition \ref{prop1} $\eta$ is a conical limit point of $X_w$. It follows by Proposition \ref{prop2} that
$\eta^{'}_1=\eta^{'}_2$. Therefore, $\eta_1$ can be flowed to $\eta_2\in \partial X_{w_2}$. Hence we are done by Theorem 
\ref{old-thm-v2}. $\Box$


\begin{example}
We now give an example contradicting the conclusion of
the Corollary $4.10$ of \cite{ps-limitset}. Suppose $\mathbb F$ is a free group on six generators, $\mathbb F=<a,b,c,x,y,z>$.
Let $H=<a,b,c>$, $K=<x,y,z>$. Suppose $\phi$ is a hyperbolic automorphism of $\mathbb F$ such that $\phi(H)=H$ and $\phi(K)=K$.
Let $G$ be the semidirect product of $F$ and $Z=<\phi>$ for the natural action of $Z$ on $\mathbb F$. Let 
$H_1=<H,\phi>$ and $K_1=<K,\phi>$. Then $H_1$, $K_1$ are both hyperbolic. (It is not difficult to see that
$H_1, K_1$ are quasiconvex subgroups of $G$.) However,  $\lim_{n\map \infty} \phi^n\in \partial G$
is a limit point of both $H$ and $K$ which are quasiconvex in $\mathbb F$ although by construction $H\cap K=(1)$.
\end{example}

{\bf Remark:} One needs to exlpore $\Lambda(G_v)\cap \Lambda(G_w)\setminus \Lambda_c(G_v\cap G_w)$.

\bibliography{lim-int.bib}
\bibliographystyle{amsalpha}

\end{document}